
\magnification=\magstep1
\pageno=0
\centerline{}
\vskip 2cm
\centerline{\bf ON THE STRUCTURE OF TENSOR PRODUCTS OF $\ell_p$-SPACES}
\vskip1cm
\centerline{by}
\vskip1cm
\centerline{Alvaro Arias${}^*$}

\centerline{The University of Texas at San Antonio}

\bigbreak
\centerline{and}
\bigbreak
\centerline{Jeff D. Farmer${}^{**}$}
\centerline{The University of Northern Colorado}
\vskip1.2in
\baselineskip18pt
\noindent{\bf Abstract.} We examine some structural properties of (injective and projective) tensor products of $\ell_p$-spaces (projections, complemented subspaces, reflexivity, isomorphisms, etc.). We combine these results with combinatorial arguments
to address the question of primarity for these spaces and their duals.
\vskip1in

\baselineskip12pt
\vfill
\item{${}^*$} This research was partially supported by NSF DMS--8921369. 

\item{${}^{**}$} Portions of this paper form a part of the second author's 
Ph.D. dissertation under the supervision of W.B. Johnson.  This research was also partially supported by NSF DMS--8921369.
\eject

\magnification=\magstep1
\pageno=1
\baselineskip18pt
\def\SN{\ell_{p_1}\hat{\otimes}\cdots \hat{\otimes}\ell_{p_N}}
\def\Sn{\ell_{p_1}\hat{\otimes}\cdots \hat{\otimes}\ell_{p_{N-1}}}
\def\Sm{\ell_{p_1}\hat{\otimes}\cdots \hat{\otimes}\ell_{p_{N-2}}}
\def\Nsigma{(\sigma_1,\sigma_2,\cdots,\sigma_N)}
\def\nsigma{(\sigma_1,\sigma_2,\cdots,\sigma_{N-1})}
\def\MN{{\bf N}^N}
\def\Mn{{\bf N}^{N-1}}
\def\N{{\bf N}}
\def\endpf{\ \ \vrule height6pt width4pt depth2pt}

\noindent{\bf 0. Introduction.}
\medbreak 
\par
A Banach space $X$ is {\it prime} if every infinite-dimensional complemented subspace contains a further subspace which is isomorphic to $X$.
A Banach space $X$ is said to be {\it primary} if whenever $X=Y\oplus Z$, 
$X$ is isomorphic to either $Y$ or $Z$.  The classical examples of prime 
spaces are the spaces $\ell_p$, $1\le p\le\infty$.  Many spaces derived 
from the $\ell_p$-spaces in various ways are primary (see for example [AEO] 
and [CL]).

\par
The primarity of $B(H)$ was shown by Blower [B] in 1990, and Arias [A] 
has recently developed further techniques which are used to prove the 
primarity of $c_1$, the space of trace class operators (this was first 
shown by Arazy [Ar1, Ar2]).  It has become clear that these techniques 
are not naturally confined to a Hilbert space context; in the present paper 
we wish to extend the results to a variety of tensor products and operator 
spaces of $\ell_p$-spaces (and in some cases 
${\cal L}_p$-spaces).  We also include some related results.
\par
Some of the intermediate propositions (on factoring operators through the 
identity) may actually be true for a wider class of Banach spaces (those with
unconditional bases which have nontrivial lower and upper estimates).  In fact, 
the combinatorial aspects of the factorization can be applied quite generally,
and may have other applications.  The proofs of primarity, however, rely on 
Pe\l czy\'nski's decomposition method which is not so readily extended.  
We have thus kept mainly to the case of injective and projective 
tensor products of $\ell_p$ spaces throughout.
The results we obtain apply to the growing study of polynomials on Banach spaces since polynomials may be considered as symmetric multilinear operators with an equivalent norm (see [FJ], [M], or [R]).
\par
Our main results are:
\medbreak
\item{(1)} If $1<p<\infty$, then $B(\ell_p)\approx B(L_p)$.
\medbreak
\item{(2)} If ${1\over p_i}+{1\over p_j}\leq1$ for every $i\neq j$, 
or if all of the $p_i$'s are equal, then $\SN$ is primary.
\medbreak
\item{(3)} $\ell_p$ embeds into $\SN$ if and only if there exists 
$A\subset \{1,2,\cdots,n\}$ such that ${1\over p}=\min\{\sum_{i\in A}{1\over p_i},1\}$.
\medbreak
\item{(4)} If $1\leq p<\infty$ and $m\geq1$, then the space of homogeneous
analytic polynomials ${\cal P}_m(\ell_p)$ and the symmetric tensor
product of $m$ copies of $\ell_p$ are primary.
\medbreak

The paper is organized as follows.  In Section 1 we set notation, definitions 
and some necessary but more or less known facts.  In Section 2 we show that 
$B(\ell_p)$, the Banach space of bounded linear operators on $\ell_p$, is 
isomorphic to $B(L_p)$, and in fact to $B(X)$ whenever $X$ is a separable 
${\cal L}_p$-space, along with some more general results we require later.  In 
Section 3 we will construct a multiplier through which
a given operator on tensor products may be factored; we then use this to show 
that some projective tensor products are primary.  In Section 4 we will prove
that the $\ell_p$ subspaces of $\SN$ are the ``obvious'' ones and use this
to prove that some projective tensor products are not primary (for example, 
$\ell_2\hat{\otimes}\ell_{1.5}$ is not primary).  Section 5
covers the question of primarity in the injective tensor products and operator
spaces, a situation not always dual to the projective case and calling for
somewhat different techniques.  Section 6 is an appendix in which we prove the 
technical lemmas we use in Section 3.
\par

We would like to thank W.B. Johnson for organizing the summer workshops in
Linear Analysis and Probability at Texas A\&M University in 1991-1993, and the
NSF for funding them.

\bigbreak
\noindent{\bf 1. Preliminaries.}
\medbreak 
Unless explicitly stated, all references to $\ell_p$ spaces will assume that
$1< p<\infty$, and will adhere the notational convention that
${1\over p_i}+{1\over q_i}=1$ or sometimes ${1\over r}+{1\over r'}=1.$

Define
$$X=\SN.$$  
We can identify its predual $X_*$ and dual $X^*$ as follows 
$$\eqalign{X_*&=\ell_{q_1}\check{\otimes}\cdots\check{\otimes}\ell_{q_N}\cr
           X^*&=B(\ell_{p_1},(\ell_{p_2}\hat{\otimes}\cdots\hat{\otimes}
                     \ell_{ p_N})^*)\cr
    &\equiv B(\ell_{p_1},B(\cdots B(\ell_{p_{N-1}},\ell_{q_N})\cdots)).\cr}$$

\noindent
The elements of $X$, $X_*$, or $X^*$ have representations as an
infinite $N$-dimensional matrix of complex numbers (we must keep in mind,
however, that this representation may not be the most efficient for
computing the tensor product norm) where the element in the 
$\alpha=(\alpha_1,\cdots,\alpha_N)\in\MN$ position is the 
coefficient of the ``matrix element'' 
$e_\alpha=e_{\alpha_1}\otimes\cdots\otimes e_{\alpha_N}$ with 
$e_{\alpha_j}$ being the $\alpha_j$-th element in the unit vector basis of $\ell_{p_j}.$
All subspaces we consider are norm-closed, and when we indicate the linear span of elements we always mean the closed span.
\par
The following elementary lemma is very important to the structure of projective
tensor products.
\medbreak
{\sl
LEMMA 1.1. Let $X$ and $Y$ be Banach spaces and $S\in B(X)$, $T\in B(Y)$.
Then $S\otimes T\in B(X\hat\otimes Y)$ is defined by
$S\otimes T(x\otimes y)=S(x)\otimes T(y)$ and satisfies
$\|S\otimes T\|\leq \|S\|\|T\|$.}
\medbreak
As a consequence of this we get that projective tensor products of 
Banach spaces with bases have bases. 
\medbreak
{\sl
PROPOSITION 1.2. Let $X$ and $Y$ be Banach spaces with bases $(e_n)_n$ and
$(f_n)_n$ respectively. Then $X\hat\otimes Y$ has a basis. Moreover, 
we take the elements of the basis from the ``shell'' 
$\partial M_n=[e_i\otimes e_j\colon \max\{i,j\}=n]$; i.e., 
$e_1\otimes f_1$,\ \  $e_2\otimes f_1,e_2\otimes f_2,e_1\otimes f_2$,\ \    
$e_3\otimes f_1,e_3\otimes f_2,e_3\otimes f_3,e_2\otimes f_3,e_1\otimes
f_3,\cdots$, etc.
}
\medbreak
The proof of this is easy. On the one hand it is clear that the span of those
vectors is dense and using Lemma 1.1 (with the operators replaced 
by projections) we see that the initial segments are uniformly complemented,
because $\partial M_n$ is clearly complemented.

As a consequence we get that $\SN$ has a basis consisting of $e_\alpha$'s.
Moreover, we can use Lemma 1.1 to prove that $\partial M_n=[e_\alpha\colon
\alpha\in\MN, \max\{\alpha_1,\cdots,\alpha_N\}=n]$ is 2-complemented and 
that $(\partial M_n)_n$ forms a
Schauder decomposition for $\SN$; we also see that $(L_\alpha)_\alpha$
is a Schauder decomposition for $\SN$ where $\alpha\in\Mn$ and
$L_\alpha=[e_\alpha\otimes e_j\colon j\in\N]$. 
(A more complete discussion of this situation appears in [R]).
We will use these facts
in Section 3.

The next theorem gives us the two most basic ingredients of our analysis.
We will prove 
that the main diagonals are 1-complemented and will identify them exactly; we will also state under what conditions the triangular 
parts of $\SN$ are
complemented. It is known that the main triangular part of 
$\ell_p\hat{\otimes}\ell_q$ is complemented if and only if
${1\over p}+{1\over q}>1$ (See [KP], [MN] and [Be]).

\medbreak
{\sl
THEOREM 1.3. Let $X=\SN$. Then the main diagonal 
${\cal D}=[e_n\otimes\cdots\otimes e_n\,:\,n\in\N]$ is 1-complemented and
satisfies ${\cal D}\equiv \ell_r$ where ${1\over r}=\min\{1,
\sum_{i=1}^N{1\over p_i}\}.$
As a consequence we get that $X\approx(\sum\oplus X)_r.$
Moreover, if $j,k$ are fixed,
then the canonical projection onto 
$[e_{i_1}\otimes e_{i_2}\otimes\cdots\otimes e_{i_N}\,:\,i_k\geq i_j]$
is bounded if and only if ${1\over p_k}+{1\over p_j}>1.$
}
\medbreak
This theorem is known for $n=2$, and in some respects for larger $n$ as well (see for example [Z]). For completeness we show here how the case $n=2$ may be extended. 
\medbreak

PROOF.  For $1<k\leq N$, let $P_{1,k}\in B(\ell_{p_1}\hat\otimes\ell_{p_k})$
be the main diagonal projection and $I_{1,k}$ be the identity on
$\ell_{p_2}\hat\otimes\cdots\hat\otimes\ell_{p_{k-1}}\hat\otimes
\ell_{p_{k+1}}\hat\otimes\cdots\hat\otimes\ell_{p_N}$. 
Then $P_{1,k}\otimes I_{1,k}$ is the projection on $\SN$ defined by
$P_{1,k}\otimes I_{1,k} e_\alpha=e_\alpha$ if $\alpha_1=\alpha_k$ and
zero otherwise.

Let $P=(P_{1,2}\otimes I_{1,2})\cdots (P_{1,N}\otimes I_{1,N})$. It is easy
to see that $Pe_\alpha=e_\alpha$ if $\alpha_1=\cdots=\alpha_N$ and
$Pe_\alpha=0$ otherwise. This tells us that ${\cal D}$ is complemented.

When $N=2$, the main diagonal of $\ell_{p_1}\hat\otimes\ell_{p_2}$
is isometric to $\ell_r$ where ${1\over r}=\min\{1,{1\over p_1}+{1\over p_2}\}$.
We apply an induction step for $N>2$. The key to the induction step is
the following: Let $D$ be the ``diagonal-projection'' on a
projective tensor products of $\ell_p$-spaces. Then it is easy to see that
$D(\SN)\equiv D( D(\Sn)\hat\otimes\ell_{p_N})$.

Notice that if the $P_{1,k}$'s above are block projections, then 
we conclude that the block
diagonal projections are bounded. By taking those to be infinite 
and using the previous paragraph, we see that 
$X\approx\bigl(\sum\oplus X\bigr)_r$.

For the last part let $T_{k,j}$ be the upper triangular projection on
$\ell_{p_k}\hat\otimes\ell_{p_j}$ and $I_{k,j}$ be the identity on
${\hat\otimes}_{i\neq k,j}\ell_{p_i}$. $T_{k,j}$ is bounded if and only if
${1\over p_k}+{1\over p_j}\leq1$. Therefore, the same is true for
$T_{k,j}\otimes I_{k,j}\in B(\SN)$. \endpf
\medbreak
REMARKS. (1) To prove that $X\approx\bigl(\sum\oplus X\bigr)_r$ we used
Pe\l czy\'nski's decomposition method. This says that if two
Banach spaces $X_1$ and $X_2$ embed complementably into each other
and if for some $1\leq p\leq\infty$, $X_1\approx\bigl(\sum\oplus X_1\bigr)_p$,
then $X_1\approx X_2$.

(2) We will work in Section 3 with 
${\cal T}=[e_\alpha\colon\alpha_1<\alpha_2<\cdots<\alpha_N]$.
Some of the results from Theorem 1.3 hold for this space. For instance,
the block projections are bounded. This implies that
${\cal T}\approx\bigl(\sum{\cal T}\bigr)_r$ where $r$ is as in Theorem 1.3.

(3) It is clear that when $r=1$ then $\SN$ is not reflexive. It is
not very difficult to prove that if $r>1$ then $\SN$ is reflexive.
\medbreak 

\vfill\eject

\bigbreak
\noindent{\bf 2. Isomorphisms of Spaces of Operators on $\ell_p$.}
\medbreak
In this section we will show that $B(\ell_p)$ is isomorphic to $B(X)$
when $X$ is any separable ${\cal L}_p$-space. In particular, $B(\ell_p)$
is isomorphic to $B(L_p[0,1])$.  A consequence of this is that $B(\ell_2)$ embeds complementably in $B(\ell_p)$ for $1<p<\infty$.
\medbreak
{\sl
THEOREM 2.1. Let $X$ and $Y$ be separable ${\cal L}_p$- and ${\cal L}_q$-spaces 
respectively with $1<p\leq q$. Then $B(X,Y)\approx B(\ell_p,\ell_q)
\approx\bigl(\sum_{n=1}^\infty\oplus B(\ell_p^n,\ell_q^n)\bigr)_\infty.$
}
\medbreak
We also obtain an isomorphic representation for
$(\ell_{p_1}\hat\otimes\cdots\hat\otimes\ell_{p_N})^*$ when 
$\sum_{i\leq N}{1\over p_i}\geq 1$.
\medbreak
{\sl
THEOREM 2.2. Let $X=\SN$ be such that 
${1\over r}=\min\{1,\sum_{i=1}^N{1\over p_i}\}=1$. 
Then $X^*\approx \bigl(\sum_{n=1}^\infty
\ell_{q_1}^n\check{\otimes}\cdots
\check{\otimes}\ell_{q_N}^n\bigr)_\infty$.
}
\medbreak
The proof of these two theorems is very similar;
they use Pe\l czy\'nski's decomposition method. 

For Theorem 2.1 notice that 
$B(\ell_p,\ell_q)\equiv (\ell_p\hat\otimes\ell_{q'})^*$ where 
${1\over q}+{1\over q'}=1$. Hence, if $p\leq q$ (i.e., 
${1\over p}+{1\over q'}\geq 1$), Theorem 1.3 tells us that
$\ell_p\hat\otimes\ell_{q'}\approx\bigl(\sum\oplus
\ell_p\hat\otimes\ell_{q'}\bigr)_1$, and then
$B(\ell_p,\ell_q)\approx\bigl(\sum\oplus B(\ell_p,\ell_q)\bigr)_\infty$.
For Theorem 2.2, notice that Theorem 1.3 implies that $X\approx\bigl(\sum
\oplus X\bigr)_1$; therefore, 
$X^*\approx\bigl(\sum\oplus X^*\bigr)_\infty$.

Then it is enough to prove that each space embeds complementably into
the other. We prove these facts for Theorem 2.1 in the next two lemmas
and indicate how to do it for Theorem 2.2 at the end of the section.

A Banach space $X$ is ${\cal L}_p$ if its finite dimensional subspaces are like
those of $\ell_p$. If $1<p<\infty$, the separable ${\cal L}_p$-spaces are
the complemented subspaces of $L_p[0,1]$ not isomorphic to $\ell_2$.

We use the following properties of a separable ${\cal L}_p$-space $X$:
(1) $X$ contains a complemented copy of $\ell_p$, and (2) There is an 
increasing (by inclusion) sequence of finite dimensional subspaces which
are uniformly isomorphic to finite dimensional $\ell_p$-spaces. Moreover,
they are uniformly complemented and their union is dense in $X$.  For more information on ${\cal L}_p$-spaces see [LP] or [JRZ].

\medbreak
{\sl
LEMMA 2.3.  Suppose that $1<p\leq q$ and let $X$ and $Y$ be
separable $ {\cal L}_p$ or ${\cal L}_q$ spaces. 
Then $B(X,Y)$ embeds
complementably in 
$W={\bigl(\sum_{n=1}^\infty \oplus B(\ell_p^n,\ell_q^n)\bigr)}_\infty.$
}
\medbreak
PROOF.  By the assumptions on $X$ and $Y$, we can find 
$\phi_n\colon B(\ell_p^n,\ell_q^n)\to B(X,Y)$ and
$\psi_n\colon B(X,Y)\to B(\ell_p^n,\ell_q^n)$ satisfying: 
(1) $\psi_n\phi_n=I_n$, the identity on $B(\ell_p^n,\ell_q^n)$, and
(2) for every $T\in B(X,Y)$, $\phi_n\psi_n(T)\to T$ in the $w^*$-topology.

Then define $\Psi:B(X,Y)\to W$ by $\Psi(T)=(\psi_n(T))_n$. Let ${\cal U}$ be
a free ultrafilter in $\N$ and define $\Phi:W\to B(X,Y)$ by
$\Phi((T_n))=\lim_{n\in{\cal U}}\phi_n(T_n)$ where the limit is taken in the
$w^*$-topology.  We can easily verify that $\Phi\Psi=I$, the identity on
$B(X,Y)$, and the conclusion follows.\endpf

\medbreak
{\sl
LEMMA 2.4. Let $X$ and $Y$ be ${\cal L}_p$ and ${\cal L}_q$-spaces 
respectively,
with $1<p\leq q$ and let $W$ be as above. 
Then $W$ embeds complementably into $B(X,Y)$
}
\medbreak
PROOF. It is clear that $W$ embeds complementably into $B(\ell_p,\ell_q)$,
because $B(\ell_p,\ell_q)$ has $\ell_\infty$-blocks down the diagonal.
Moreover, if $X$ is a separable ${\cal L}_p$-space, $1<p<\infty$, then
$X$ contains a complemented copy of $\ell_p$. Since the same is true
for $Y$ we see that $B(\ell_p,\ell_q)$ embeds complementably into
$B(X,Y)$.\endpf
\medbreak
REMARK. For Theorem 2.2 notice that
$Z_n=B(\ell_{p_1}^n,B(\cdots B(\ell_{p_{N-1}}^n,\ell_{q_N}^n)\cdots))$ 
is isometric to
$\ell_{q_1}^n\check\otimes\cdots\check\otimes\ell_{q_N}^n$ and is
1-complemented in $X^*$. 
(We use this to show that $\bigl(\sum_n\oplus Z_n\bigr)_\infty$ embeds
complementably into $X^*$).
Moreover, $\bigcup_nZ_n$ is $w^*$-dense in $X^*$. (We may use this and
an ultrafilter argument to show the reverse complemented inclusion).
\medbreak

\bigbreak
\noindent{\bf 3. Primarity of Projective Tensor Products.}
\medbreak

We devote most of this section to the proof of the following theorem.

\medbreak
{\sl
THEOREM 3.1. Let $X=\SN$ be such that  
${1\over p_i}+{1\over p_j}\leq1$ for every $i\neq j$.
Then $X$ is primary.
}
\medbreak
The proof of this theorem will follow easily from the next proposition that
was inspired by results of Blower [B] and was used in [A] 
in a similar context.
The ideas involved in this ``factorization'' approach are well-known (see for example Bourgain [Bo]).

We have to introduce some notation.

Let $X=\SN$, $\alpha=(\alpha_1,\cdots,\alpha_N)\in\MN$ and denote by
$e_\alpha=e_{\alpha_1}\otimes e_{\alpha_2}\otimes\cdots\otimes e_{\alpha_N}.$
Then $X=[e_\alpha\,:\,\alpha\in\MN].$  We also define 
$|\alpha|=\max\{\alpha_1,\cdots\alpha_2\};$
and introduce an order between different multiindices. Let
$\alpha\in{\bf N}^k$ and $\beta\in{\bf N}^m$; we say that 
$\alpha<\beta$ if $\max\{\alpha_1,\cdots,\alpha_k\}<
\min\{\beta_1,\cdots,\beta_m\}.$

Let $\sigma_1,\sigma_2,\cdots\sigma_N:\N\to\N$ be increasing functions 
(it will also be useful to think of the $\sigma_i$'s as infinite subsets of $\N$);
and let $\sigma=\Nsigma$ be a function on $\MN$ defined by
$\sigma(\alpha)=(\sigma_1(\alpha_1),\cdots,\sigma_N(\alpha_N)).$
Then define 
$$\eqalign{J_\sigma&:\SN\to\SN\cr
           K_\sigma&:\SN\to\SN\cr}$$
by 
$$J_\sigma e_\alpha=e_{\sigma(\alpha)}\quad\hbox{and}\quad
K_\sigma e_\alpha=\cases{e_\beta &if there exists $\beta$ such that 
                                   $\sigma(\beta)=\alpha$\cr
                         0 &otherwise.}$$

$J_\sigma$ and $K_\sigma$ have many important algebraic properties. 
$J_\sigma$ is
one-to one, $K_\sigma$ is onto and 
$K_\sigma   J_\sigma=I.$
Moreover, they compose nicely; that is, if $\sigma=\Nsigma$ and 
$\psi=(\psi_1,\cdots,\psi_N)$, then
$$J_\sigma  J_\psi=J_{\sigma \psi}\quad\hbox{ and }\quad
 K_\sigma  K_\psi=K_{\sigma \psi}.$$

We are now ready to state the proposition.
\medbreak
{\sl
PROPOSITION 3.2. Let ${1\over p_i}+{1\over p_j}\leq1$ for every $i\neq j$. Then if $\Phi\in B(\SN)$ and $\epsilon>0$
there exist $\sigma=\Nsigma$ and $\lambda\in{\bf C}$ such that
$$\|K_\sigma\Phi J_\sigma-\lambda I\|<\epsilon.$$
Thus one of $K_\sigma \Phi J_\sigma$ or  $K_\sigma(\Phi-I) J_\sigma$
is invertible.
}
\medbreak
It is immediate from this proposition that $I$, the identity on $\SN$,
factors through $\Phi$ or through $I-\Phi$ which implies trivially
that if $\SN\approx X\oplus Y$ then $\SN$ embeds complementably
into $X$ or $Y$. Since $\SN$ is isomorphic to its infinite $r$-sum,
the Pe\l czy\'nski
decomposition method implies that $\SN$ is primary.
\medbreak

We will present a sketch of the proof. For $\Phi\in\SN$ and $\alpha\in\MN$ we have, 
$$\Phi e_\alpha=\sum_{\beta\in\MN}\lambda_{\alpha,\beta}e_\beta,$$
for some $\lambda_{\alpha,\beta}\in{\bf C}$. 
Our goal will be to come with a series of the aforementioned $J$-maps and $K$-maps which will allow us to get $K\Phi J\approx\lambda I$. We
will do this is several steps, fixing progressively more restrictive
portions of the range of $\beta$. We can do this since this maps compose nicely; however we must be careful no to destroy previous work (see the assumption below).
More precisely, Step 1 asserts that we can find $K_1, J_1$ such that 
$K_1\Phi J_1\approx\Phi_1$ and for every $n\in\N$,
$$\alpha\in\MN,\quad |\alpha|=n\quad\quad\Longrightarrow\quad\quad \Phi_1 e_\alpha=\sum_{|\beta|=n}
\lambda_{\alpha,\beta}^{(1)}e_\beta,$$
for some $\lambda_{\alpha,\beta}^{(1)}\in{\bf C}$. This is clearly an improvement in the range of $\beta$, but we still have that
$\{\beta\in\MN : |\beta|=n\}$ is a big set. After Steps 2, 3 and 4 we have 
$K_4, J_4$ such that $K_4\Phi J_4\approx \Phi_4$ and for every $n\in\N$, $j\in\N$,
$$
\alpha\in\Mn,\quad |\alpha|=n,\quad |\alpha|<j
\quad\quad
\Longrightarrow
\quad\quad 
\Phi_4 (e_\alpha\otimes e_j)=\sum_{|\beta|=n\atop \beta\in\Mn}
\lambda_{\alpha,\beta}^{(4)}e_\beta\otimes e_j.
$$
Step 5 gives us $K_5, J_5$ such that
$K_5\Phi J_5\approx\Phi_5$ and for every $n\in\N$, $\gamma\in\N^2$,
$$
\alpha\in\N^{N-2},\quad |\alpha|=n,\quad \alpha<\gamma
\quad\quad
\Longrightarrow
\quad\quad 
\Phi_5 (e_\alpha\otimes e_\gamma)=\sum_{|\beta|=n\atop \beta\in\N^{N-2}}
\lambda_{\alpha,\beta}^{(5)}e_\beta\otimes e_\gamma.
$$
Finally Step 6 provides the general induction argument.
\medbreak
We will apply our arguments on
${\cal T}=[e_\alpha\,:\,
\alpha_1<\alpha_2<\cdots<\alpha_N]$ without loss of generality in order to simplify notation, keeping in mind that they will be repeated many times when
the order of the $\alpha_i$'s is different. We will choose $\sigma$ so that
$J_\sigma$ and $K_\sigma$ ``respect'' that order. More precisely, consider the permutation group $\Pi_n$ and a multiindex 
$\alpha=(\alpha_1,\cdots,\alpha_N)\in\MN$.  
We choose $\sigma$ so that the (not necessarily complemented) subspaces
${\cal T}(\pi)=
[e_\alpha\,:\,\alpha_{\pi(1)}<\alpha_{\pi(2)}<\dots<\alpha_{\pi(N)}]$
are invariant for $J_\sigma$ and $K_\sigma$; i.e., 
$J_\sigma {\cal T}(\pi)\subset {\cal T}(\pi)$ and
$K_\sigma {\cal T}(\pi)\subset {\cal T}(\pi)$.
Notice that the ${\cal T}(\pi)$'s 
``exhaust''  the $N$-dimensional matrix array on which we represent $\SN$
(modulo diagonal elements, which we always ignore; see Step 2).
\medbreak
{\sl
ASSUMPTION. Assume from now on that whenever we choose $\sigma=\Nsigma$, it always ``preserves the order'',
that is,
if $i<j,$ then $\sigma_k(i)<\sigma_l(j),$
for every $k,l\leq N$.}

We can always satisfy this assumption by passing to subsequences whenever we are choosing the sets $\sigma$, which our technical lemmas allow us to do.
\medbreak
EXAMPLE: It might be instructive to consider the following example ``far''
from a multiplier. Let $\Phi:\ell_2\hat{\otimes}\ell_2\to \ell_2\hat{\otimes}\ell_2$ be the transpose operator, i.e.,
$\Phi e_i\otimes e_j=e_j\otimes e_i$.
Then choose $\sigma_1$ the set of even integers,
$\sigma_2$ the set of
odd integers and $\sigma=(\sigma_1,\sigma_2)$. We verify easily that $K_\sigma\Phi J_\sigma=0$ thus satisfying the conclusion of Proposition 3.2.
\medbreak

Our steps require the repeated use of two technical lemmas whose proof we delay until Section 6.

\bigbreak

{\sl
STEP 1. Let $\partial M_n=[e_\alpha\,:\,\alpha\in \MN,
|\alpha|=n]$ with projection $Q_n$. Then for every 
$\Phi\in B(\SN)$ and $\epsilon>0$, there
exist $\sigma=\Nsigma$ and $\Phi_1\in B(\SN)$ such that
$\|\Phi_1-K_\sigma\Phi J_\sigma\|<\epsilon$ and for every $n$,
$\Phi_1\partial M_n\subset \partial M_n$.
}
\medbreak

The proof of this step is an immediate consequence of the following lemma.
Remember that $(\partial M_n)_n$ forms
a Schauder decomposition for $\SN$.

\medbreak
{\sl
BASIC LEMMA 1. Let $\Phi\in B(\SN)$.  
Then for every $\epsilon_{n,m}>0$ we can find
$\sigma=\Nsigma$ such that if
$x\in\partial M_n$, and $n\neq m$, then
$\|Q_m K_\sigma\Phi J_\sigma x\|\leq\epsilon_{n,m}\|x\|.$
}
\medbreak

We prove Basic Lemma 1 in the appendix (if $X=\ell_p$ the proof is very easy).

Choose $\epsilon_{n,m}$ in Basic Lemma 1 so that $\epsilon_n=\sum_{m=1}^\infty\epsilon_{n,m}$ and $\sum_{n=1}^\infty\epsilon_n<{\epsilon\over2}$. Then
define $\Phi_1$ on $\SN$ as follows: For $x\in \partial M_n$, let
$$\Phi_1(x)=Q_nK_\sigma\Phi J_\sigma x.$$
If $x\in\partial M_n$, then $\|(\Phi_1-K_\sigma\Phi J_\sigma)x\|\leq\epsilon_n\|x\|$. If 
$x\in\SN$, we have that $x=\sum_{n=1}^\infty x_n$ where $x_n\in\partial M_n$ and $\|x_n\|\leq2\|x\|$. Therefore,
$$\|(\Phi_1-K_\sigma\Phi J_\sigma)x\|\leq\sum_{n=1}^\infty
\|(\Phi_1-K_\sigma\Phi J_\sigma)x_n\|<\epsilon\|x\|.\endpf$$

\medbreak

{\sl
STEP 2. Let $\Phi\in B(\SN)$ be such that 
$\Phi\partial M_n\subset\partial M_n$ for every $n$, then we can find 
$\sigma=\Nsigma$ such that $\Phi_2=K_\sigma\Phi J_\sigma$ ``respects''
the place where $\alpha\in\MN$ takes is maximum; for example, if
the maximum takes place in the last coordinate, i.e., 
$\alpha\in\Mn$ and $|\alpha|<j$, then
$$\Phi_2 (e_\alpha\otimes e_j)=\sum_{|\beta|<j}
\lambda_{\alpha,\beta,j}e_\beta\otimes e_j,$$
and we also have similar results for the other coordinates.
}
\medbreak

We attain this by ``disjointifying'' the different faces.
For $i\leq N$ let $\sigma_i(j)=N(j-1)+i$, and $\sigma=\Nsigma$.
It is easy to see that $\Phi_2=K_\sigma\Phi J_\sigma$ satisfies the
required property.
Indeed, if $\alpha\in\Mn$ and
$|\alpha|<n$, then $e_\alpha\otimes e_n\in\partial M_n$, and $J_\sigma (e_\alpha\otimes e_n)\in\partial M_{\sigma_N(n)}$. Hence,
$$\Phi J_\sigma (e_\alpha\otimes e_n)=\sum_{|\gamma|=\sigma_N(n)}
\lambda_{\sigma(\alpha,n),\gamma} e_\gamma.$$
Recall that $K_\sigma e_\gamma=e_\eta$ if $\sigma(\eta)=\gamma$ for some $\eta$
and $K_\sigma e_\gamma=0$ otherwise. Since
the ranges of the $\sigma_i$'s are disjoint,
$\sigma_N(n)$ is nonzero only for the last coordinate. Therefore, if 
$|\gamma|=\sigma_N(n)$ and $\sigma(\eta)=\gamma$, the last coordinate of $\eta$
must be $n$; i.e., $e_\eta=e_\beta\otimes e_n$ for some $\beta\in\Mn$,
and since $\sigma$ preserves the order, $|\beta|<n$. That is,
$$K_\sigma\Phi J_\sigma (e_\alpha\otimes e_n)=\sum_{|\beta|<n}
\lambda_{\sigma(\alpha,n),\sigma(\beta,n)}e_\beta\otimes e_n.$$
We denote $\lambda_{\sigma(\alpha,n),\sigma(\beta,n)}$ by $\lambda_{\alpha,\beta,n}$.\endpf

\medbreak
To make the notation a bit clearer we will state the hypothesis 
and the conclusion of the
steps when the maximum takes place in the $N$th coordinate. However
the other cases are identical and we will assume that (after repeating
the step for the other coordinates) the same result holds for these cases.
\medbreak
{\sl
STEP 3. Let $\Phi\in B(\SN)$ be such that whenever $\alpha\in\Mn$, $j\in\N$ satisfy $|\alpha|<j$, then
$\Phi (e_\alpha\otimes e_j)=\sum_{|\beta|<j}
\lambda_{\alpha,\beta,j}e_\beta\otimes e_j$.
Then for every $\epsilon>0$ there exist $\sigma=\Nsigma$ and 
$\Phi_3\in B(\SN)$ such that $\|\Phi_3-K_\sigma\Phi J_\sigma\|<\epsilon$ and
whenever $\alpha\in\Mn$, $j\in\N$ satisfy $|\alpha|<j$, then
$$\Phi_3(e_\alpha\otimes j)=\sum_{|\beta|<j}
\mu_{\alpha,\beta}e_\beta\otimes j.$$
}
\medbreak
The proof of this step follows from the next lemma.
\medbreak
{\sl
BASIC LEMMA 2. Let $\Phi\in B(\SN)$ be such that whenever $\alpha\in\Mn$, $j\in\N$ satisfy $|\alpha|<j$, then 
$\Phi(e_\alpha\otimes e_j)=\sum_{|\beta|<j}\lambda_{\alpha,\beta,j}
e_\beta\otimes e_j.$
Then for every $\epsilon_{\alpha,\beta,j}>0$ with $j>\max\{|\alpha|,|\beta|\}$,
we can find $\sigma=\Nsigma$ (respecting the order) such that if we set
$\tilde{\sigma}=\nsigma$ then
$$\eqalign{&\lim_{j\to\infty}\lambda_{\tilde{\sigma}(\alpha),
\tilde{\sigma}(\beta),\sigma_N(j)}=\lambda_{\tilde{\sigma}(\alpha),
\tilde{\sigma}(\beta)};\cr
&|\lambda_{\tilde{\sigma}(\alpha),
\tilde{\sigma}(\beta),\tilde{\sigma}_N(j)}-
\lambda_{\tilde{\sigma}(\alpha),
\tilde{\sigma}(\beta)}|\leq\epsilon_{\alpha,\beta,j}.\cr}$$
}
\medbreak

We also give the proof of Basic Lemma 2 in the appendix.
Then set $\tilde{\Phi}=K_\sigma \Phi J_\sigma$, and let 
$L_\alpha=[e_\alpha\otimes e_j\,:\,j\in {\bf N}]$ with projection 
$P_\alpha$. Since
$\Sn$ has a basis consisting of $e_\alpha$'s, we have that
$(L_\alpha)_\alpha$ forms a Schauder decomposition for $\SN$. 

Define $\Phi_3\in B(\SN)$ by
$$P_\beta\Phi_3(e_\alpha\otimes e_j)=
\cases{\lambda_{\tilde{\sigma}(\alpha),\tilde{\sigma}(\beta)}
                                             e_\beta\otimes e_j& if 
                          $\max\{|\alpha|,|\beta|\}<j;$\cr
              P_\beta\tilde{\Phi} (e_\alpha\otimes e_j)& otherwise.\cr}$$
      
Let $\alpha,\beta\in\MN$;
$\epsilon_{\alpha,\beta}=\sum_{j>\max\{|\alpha|,|\beta|\}}
\epsilon_{\alpha,\beta,j}$;
and $x\in L_\alpha$; i.e., 
$x=\sum_{j=1}^\infty e_\alpha\otimes c_je_j$. Then,
$$P_\beta\Phi_3 x-P_\beta\tilde{\Phi}x=
\sum_{{j>\max\{|\alpha|,|\beta|\}}}
e_\beta\otimes (\lambda_{\tilde{\sigma}(\alpha),\tilde{\sigma}(\beta)}
-\lambda_{\tilde{\sigma}(\alpha),\tilde{\sigma}(\beta),\sigma_N(j)})
c_j e_j.$$
Hence,
$$\|P_\beta\Phi_3x-P_\beta\tilde{\Phi}x\|\leq
\sum_{j>\max\{|\alpha|,|\beta|\}}\epsilon_{\alpha,\beta,j}\max|c_j|\leq
\epsilon_{\alpha,\beta}\|x\|.$$
If we choose $\sum_\alpha\sum_\beta \epsilon_{\alpha,\beta}<\epsilon$
small enough, $\Phi_3$ is well defined and satisfies
the required properties.\endpf
\medbreak
{\sl
STEP 4.  Let $\Phi\in B(\SN)$ be such that whenever $\alpha\in\Mn$, $j\in\N$ satisfy $|\alpha|<j$, then
$\Phi (e_\alpha\otimes e_j)=\sum_{|\beta|<j}\lambda_{\alpha,\beta}
e_\beta\otimes e_j.$ Then for every
$\epsilon>0$, there exist $\sigma=\Nsigma$ and $\Phi_4\in B(\SN)$ such that
$\|\Phi_4-K_\sigma\Phi J_\sigma\|<\epsilon$ and whenever $\alpha\in\Mn$, $j\in\N$ satisfy $|\alpha|<j$, we have
$$\Phi_4 (e_\alpha\otimes e_j)=\sum_{|\beta|=|\alpha|}\mu_{\alpha,\beta}e_\beta
\otimes e_j.$$
}
\medbreak

Define $\Psi\in B(\Sn)$ by
$$(\Psi e_\alpha ,e_\beta)=\lambda_{\alpha,\beta}.$$
Since $(\Psi e_\alpha,e_\beta)=\lim_{j\to\infty}
(\Phi e_\alpha\otimes e_j,e_\beta\otimes e_j)$, $\Psi$ is a bounded map.

\medbreak
REMARK. Ideally we would like to apply an induction step and
replace $\Psi$, after a factorization of the form $K\Psi J$, by a 
multiple of the identity and then combine this with $\Phi$.  Controlling the norm of the perturbation requires a more delicate argument, however.
\medbreak
Apply Basic Lemma 1 to $\Psi$ with its respective $\partial M_n$ and
projections $Q_n$; then find $\sigma=\nsigma$ such that whenever 
$x\in\partial M_n$, and $m\neq n$,
$$\|Q_m K_\sigma\Psi J_\sigma x\|\leq \epsilon_{n,m}\|x\|.$$ 
Let $\tilde{\Psi}=K_\sigma\Psi J_\sigma$,
$\tilde{\sigma}=(\sigma_1,\cdots,\sigma_{N-1},\sigma_{N-1})$ and
$\tilde{\Phi}=K_{\tilde{\sigma}}\Phi J_{\tilde{\sigma}}$. 
Denote by
$P_n=Q_1+\cdots +Q_n$ the projection onto $[e_\alpha\,:\,
\alpha\in\Mn, |\alpha|\leq n]$. Notice now that if $|\alpha|<j$ then 
$\tilde{\Phi} (e_\alpha\otimes e_j)=(P_j \tilde{\Psi} e_\alpha)\otimes
e_j$.

Let $L_\alpha=[e_\alpha\otimes e_j\,:\,j\in{\bf N}]$. Then as we explained
after the Basic Lemma 2, $(L_\alpha)_\alpha$ forms a Schauder decomposition
for $\SN$. 

Define $\Phi_4$ on $\SN$ as follows
$$\Phi_4 (e_\alpha\otimes e_j)=
\cases{\tilde{\Phi}(e_\alpha\otimes e_j) &if $|\alpha|\geq j$;\cr
       (Q_{|\alpha|}\tilde{\Psi} e_\alpha)\otimes e_j &if $|\alpha|<j.$\cr}$$

Let $n\in\N$; $\alpha\in\Mn$ with $|\alpha|=n$; 
$\epsilon_n=\sum_{m=1,m\neq n}^\infty\epsilon_{n,m}$ and $x\in L_\alpha$;
i.e., $x=\sum_{j=1}^\infty e_\alpha\otimes c_j e_j$.
Then,
$$\eqalign{\tilde{\Phi}(x)-\Phi_4(x)
&=\sum_{j>n}(P_j\tilde{\Psi} e_\alpha-Q_n\tilde{\Psi}e_\alpha)\otimes c_je_j\cr
&=\sum_{j>n}\sum_{k=1\atop k\neq n}^j (Q_k\tilde{\Psi} e_\alpha)\otimes c_j e_j\cr
&=\sum_{k=1\atop k\neq n}^\infty (Q_k\tilde{\Psi} e_\alpha)\otimes
\sum_{j>\max\{k,n\}}c_j e_j.\cr}$$
 
Since
$\|\sum_{j>\max\{k,n\}}c_j e_j\|\leq \|x\|$, we have that
$$\|\Phi_4(x)-\tilde{\Phi}(x)\|\leq \epsilon_n\|x\|.$$

Since card$\{\alpha\,:\,|\alpha|=n\}$ is finite, it is enough to choose 
$\epsilon_n$ so that 
$\sum_{n=1}^\infty\hbox{ card}\,\{\alpha\,:\,|\alpha|=n\}\epsilon_n<\epsilon$
to insure that $\Phi_4$ is well defined and satisfies 
the required properties.\endpf
\medbreak

{\sl STEP 5. Let $\Phi\in B(\SN)$ be such that whenever $\alpha\in\Mn$ and
$j\in\N$ satisfy $|\alpha|<j$, then
$\Phi (e_\alpha\otimes e_j)=\sum_{|\beta|=|\alpha|}\lambda_{\alpha,\beta}
e_\beta\otimes e_j$. Then for every $\epsilon>0$, there exist
$\sigma=\Nsigma$ and $\Phi_5\in B(\SN)$ such that 
$\|\Phi_5-K_\sigma\Phi J_\sigma\|<\epsilon$ and whenever $i,j\in\N$ and
$\alpha\in\N^{N-2}$ satisfy $|\alpha|<i,|\alpha|<j$, then
$$\Phi_5 (e_\alpha\otimes e_i\otimes e_j)=\sum_{|\beta|=|\alpha|}
\mu_{\alpha,\beta} e_\beta\otimes e_i\otimes e_j.$$}

\medbreak
PROOF. Disjointifying for $\Sn$ as in Step 2, we can assume without
loss of generality that whenever $i,j\in\N$, $\alpha\in\N^{N-2}$
satisfy $|\alpha|<i<j$, then
$$\Phi (e_\alpha\otimes e_i\otimes e_j)=
\sum_{|\beta|<i }\lambda_{\alpha,\beta,i}
e_\beta\otimes e_i\otimes e_j.$$

Apply Basic Lemma 2 to the sequence $\{\lambda_{\alpha,\beta,i}\}$
and assume that (after factoring $\Phi$ through $K_\sigma\Phi J_\sigma$ and
renaming it $\Phi$ again) this sequence satisfies the conclusions of that lemma.

Let $L_\alpha=[e_\alpha\otimes e_i\otimes e_j\,:\,i,j\in{\bf N}]$ with
projection $P_\alpha$. Since $\Sm$ has a basis consisting of $e_\alpha$'s,
then $(L_\alpha)$ forms a Schauder decomposition for $\SN$.

Define $\tilde{\Phi}\in B(\SN)$ as follows:
$$P_\beta\tilde{\Phi}(e_\alpha\otimes e_i\otimes e_j)=
\cases{
   \lambda_{\alpha,\beta}e_\beta\otimes e_i\otimes e_j 
       &if $|\alpha|\vee|\beta|<i<j$;\cr
   P_\beta\Phi (e_\alpha\otimes e_i\otimes e_j)
       &otherwise.\cr}$$

Let $\alpha,\beta\in{\bf N}^{N-2}$; 
$\epsilon_{\alpha,\beta}=\sum_{i>|\alpha|\vee|\beta|}\epsilon_{\alpha,\beta,i}$  and 
$x\in L_\alpha$; i.e., $x=\sum_i\sum_j c_{i,j} e_\alpha\otimes e_i
\otimes e_j$. Then,
$$\eqalign{P_\beta\Phi x-P_\beta\tilde{\Phi}x
&=\sum_{i>|\alpha|\vee|\beta|}\biggl(\sum_{j=i+1}^\infty c_{i,j}[P_\beta\Phi (e_\alpha\otimes e_i\otimes
                     e_j)-P_\beta\tilde{\Phi}(e_\alpha\otimes e_i\otimes e_j)]\biggr)\cr
&=\sum_{i>|\alpha|\vee|\beta|}(\lambda_{\alpha,\beta,i}-\lambda_{\alpha,\beta})e_\beta
    \otimes e_i\otimes\biggl( \sum_{j=i+1}^\infty c_{i,j} e_j\biggr).\cr}$$
Since $\|\sum_{j=i+1}^\infty c_{i,j} e_j\|\leq\|x\|$,
$$\|P_\beta\Phi x-P_\beta\tilde{\Phi}x\|\leq\epsilon_{\alpha,\beta}\|x\|.$$

If we choose $\epsilon_{\alpha,\beta}$ small enough so that 
$\sum_{\alpha,\beta}\epsilon_{\alpha,\beta}<{\epsilon\over 2}$, then $\tilde{\Phi}$ is well defined and satisfies 
$\|\Phi-\tilde{\Phi}\|<{\epsilon\over2}$; moreover, whenever $|\alpha|<i<j$, we have
$$\tilde{\Phi}(e_\alpha\otimes e_i\otimes e_j)=
\sum_{|\beta|<i}\lambda_{\alpha,\beta}e_\beta\otimes e_i\otimes e_j.$$

Let $T:L_\alpha\to L_\beta$ be defined by $T(x)=P_\beta\tilde{\Phi}x$.
$T$ is clearly a bounded map and 
$L_\alpha\equiv L_\beta\equiv \ell_{p_{N-1}}\hat{\otimes}\ell_{p_N}$.
It follows that if $|\alpha|\vee|\beta|<i<j$ then $T(e_\alpha\otimes e_i\otimes e_j)=
\lambda_{\alpha,\beta}e_\beta\otimes e_i\otimes e_j$. Since all
the arguments work if the maximum is attained at the $(N-1)$-st coordinate
and the next maximum is attained in the last coordinate, we can also
assume that if $|\alpha|\vee|\beta|<j<i$ then
$T(e_\alpha\otimes e_i\otimes e_j)=
\mu_{\alpha,\beta}e_\beta\otimes e_i\otimes e_j$.
Thus $T$ takes value $\lambda_{\alpha,\beta}$ in the upper triangular
part of a copy of $\ell_{p_{N-1}}\hat{\otimes}\ell_{p_N}$ and the value
$\mu_{\alpha,\beta}$ in the lower part. Since we assumed that
$${1\over p_{N-1}}+{1\over p_N}\leq 1,$$
we have that $\lambda_{\alpha,\beta}=\mu_{\alpha,\beta}.$
(If $\lambda_{\alpha,\beta}\neq\mu_{\alpha,\beta}$, then
$(T-\mu_{\alpha,\beta}I)/(\lambda_{\alpha,\beta}-\mu_{\alpha,\beta})$ would 
be a projection onto the upper triangular part of 
$\ell_{p_{N-1}}\hat\otimes\ell_{p_N}$, contradicting Theorem 1.3).

Let $\alpha=(\alpha_1,\cdots,\alpha_{N-2})$ and
$\gamma=(\alpha_{N-1},\alpha_N)$. We now have that if $\alpha<\gamma$, then
$$\tilde{\Phi}(e_\alpha\otimes e_\gamma)=
\sum_{\beta<\gamma}\lambda_{\alpha,\beta} e_\beta\otimes e_\gamma.$$
Define a map $\Psi\in B(\Sm)$, as in Step 4, by
$$(\Psi e_\alpha,e_\beta)=\lambda_{\alpha,\beta}.$$
Since $\Psi$ is bounded, we apply Basic Lemma 1 to it and assume
without loss of generality (after factoring $K_\sigma\tilde{\Phi} J_\sigma$
and then renaming it $\tilde{\Phi}$ again) that if $x\in\partial M_n$ (here 
$\partial M_n$ is a subset of $\Sm$) and $m\not=n$,  
$$\|Q_m\Psi x\|\leq\epsilon_{n,m}\|x\|.$$

Let $L_\alpha=[e_\alpha\otimes e_i\otimes e_j\,:\,i,j\in{\bf N}]$ with
projection $P_\alpha$, and define $\Phi_1\in B(\SN)$ as follows
$$\Phi_5 (e_\alpha\otimes e_\gamma)=
\cases{(Q_{|\alpha|}\Psi e_\alpha)\otimes e_\gamma &if $\alpha<\gamma$\cr
        \tilde{\Phi}(e_\alpha\otimes e_\gamma) &otherwise.\cr}$$

Let $x\in L_\alpha$, $|\alpha|=n$; i.e.,  $x=\sum_\gamma c_\gamma e_\alpha\otimes e_\gamma$.
Hence,
$$\eqalign{\tilde{\Phi}x-\Phi_5x
&=\sum_{\gamma>\alpha}[c_\gamma\tilde{\Phi}e_\alpha\otimes e_\gamma-
                       c_\gamma(Q_n\Psi e_\alpha)\otimes e_\gamma]\cr
&=\sum_{k=1\atop k\neq n}^\infty(Q_k\Psi e_\alpha)\otimes
       \biggl(\sum_{\gamma>n,k}c_\gamma e_\alpha\otimes e_\gamma\biggr).\cr}$$
Since $\|\sum_{\gamma>n,k}c_\gamma e_\alpha\otimes e_\gamma\|\leq \|x\|$, the
result follows.\endpf
\medbreak
The induction step is an extension of Step 5.
\medbreak
{\sl
STEP 6. Let $\Phi\in B(\SN)$ be such that whenever 
$\alpha\in\N^k$, $\gamma\in\N^{N-k}$ satisfy $\alpha<\gamma$, then
$\Phi (e_\alpha\otimes e_\gamma)=\sum_{|\beta|=|\alpha|}\lambda_{\alpha,\beta}
e_\beta\otimes e_\gamma$. Then for every $\epsilon>0$, there exist
$\sigma=\Nsigma$ and $\Phi_k\in B(\SN)$ such that 
$\|\Phi_k-K_\sigma\Phi J_\sigma\|<\epsilon$ and whenever 
$\alpha\in{\bf N}^{k-1}$ and $\gamma\in\N^{N-k+1}$ satisfy
$\alpha<\gamma$, we have 
$$\Phi_k (e_\alpha\otimes e_\gamma)=\sum_{|\beta|=|\alpha|}
\mu_{\alpha,\beta} e_\beta\otimes e_\gamma.$$
}

SKETCH OF PROOF. Disjointifying as in Step 2 we assume that whenever
$\alpha\in{\bf N}^{k-1}$, $i\in\N$ and $\gamma\in{\bf N}^{N-k}$ satisfy
$\alpha<i<\gamma$, then
$$\Phi (e_\alpha\otimes e_i\otimes e_\gamma)=\sum_{\beta<i}
\lambda_{\alpha,\beta,i}e_\beta\otimes e_i \otimes e_\gamma.$$
Assume also that the sequence $\{\lambda_{\alpha,\beta,i}\}$ satisfies the conclusion of Basic Lemma 2. 

Let $\alpha\in\N^{k-1}$; $L_\alpha=[e_\alpha\otimes e_i\otimes
e_\gamma : i\in{\bf N}, \gamma\in\N^{N-k}]$ with projection $P_\alpha$ and
define $\tilde{\Phi}$ as in Step 5.  Since for every $i_0$,
$\|\sum_{\gamma>i_0}c_{i_0,\gamma}e_{i_0}\otimes e_\gamma\|\leq
\|\sum_{i,\gamma}c_{i,\gamma}e_i\otimes e_\gamma\|$, then
$\|\Phi-\tilde{\Phi}\|<\epsilon$ and 
whenever $\alpha\in \N^{k-1}$, $\gamma\in \N^{N-k}$ satisfy
$\alpha<i<\gamma$, then
$$\tilde{\Phi}(e_\alpha\otimes e_i\otimes e_\gamma)=
\sum_{\beta<i}\lambda_{\alpha,\beta}e_\alpha\otimes e_i\otimes e_\gamma.$$

Fix $\alpha,\beta\in\N^{k-1}$ and define $T:L_\alpha\to L_\beta$ by
$Tx=P_\beta\tilde{\Phi}x$. Since $T$ is bounded and 
$L_\alpha\equiv L_\beta\equiv \ell_{p_k}\hat\otimes\cdots\hat\otimes\ell_{p_N}$
we assume that $T$ is defined on 
$Z=\ell_{p_k}\hat\otimes\cdots\hat\otimes\ell_{p_N}$.

Decompose $Z$ into $(E_j)_{j=k}^N$, where 
$E_j=[e_\theta:\theta_j<\theta_i\hbox{ for every } i\neq j]$;
(i.e., $E_j$ is the span of those $e_\theta$ where the minimum occurs
at the $j$th coordinate).
For instance, if $e_\theta\in E_k$, then $e_\theta=e_i\otimes e_\gamma$
for some $i<\gamma$ and hence
$Te_\theta=\lambda_{\alpha,\beta}e_\theta=\lambda^{(k)}e_\theta$. 
Since all the arguments work for the other permutations of the 
coordinates we can assume
that there exist $\lambda^{(j)}$ such that if $x\in E_j$, then
$Tx=\lambda^{(j)}x$, $j=k,k+1,\cdots,N$. 

We will use that ${1\over p_i}+{1\over p_j}\leq1$ for every $i\neq j$
to conclude that the $\lambda^{(j)}$'s have to be equal. Indeed, let
$\bar{m}=(m,m,\cdots,m)\in\N^{N-k-3}$ and consider
$K_m=[e_i\otimes e_j\otimes e_{\bar{m}}:i,j\leq m]$. It is clear that
$K_m\equiv\ell_{p_k}^m\hat\otimes\ell_{p_{k+1}}^m$ and that $T$ restricted to
it gives us $\lambda^{(k)}$ in the upper triangular part and
$\lambda^{(k+1)}$ in the lower one. If $\lambda^{(k)}\neq\lambda^{(k+1)}$,
we would have that the $m$-triangular parts are uniformly
complemented and this is not true.  A similar
argument proves that the $\lambda^{(j)}$'s are all equal. 

In conclusion, if $\alpha\in\N^{k-1}$, $\gamma\in\N^{N-k+1}$,
and $\alpha<\gamma$, then
$$\tilde{\Phi}(e_\alpha\otimes e_\gamma)=\sum_{\beta<\gamma}\lambda_{\alpha,\beta} e_\beta\otimes
e_\gamma.$$

Define $\Psi\in B(\ell_{p_1}\hat\otimes\cdots\hat\otimes\ell_{p_{k-1}})$ by
$(\Psi e_\alpha,e_\beta)=\lambda_{\alpha,\beta}$; apply Basic Lemma 1
to it and finish the proof as in Step 5.\endpf
\medbreak

Iterating Step 6 we finish the proof of the proposition.
\medbreak

We will see in the next section that, for most cases, if 
${1\over p_i}+{1\over p_j}>1$, $X$ is not primary. This is not always true,
however.

\medbreak
{\sl
THEOREM 3.3. Let $1\leq p<\infty$ and $n\in\N$. Then $X=\ell_p\hat\otimes\cdots
\hat\otimes\ell_p$ ($n$ times) is primary.
}
\medbreak
PROOF. We divide the proof into two cases. If ${2\over p}\leq 1$,
this is a particular case of Theorem 3.1. If ${2\over p}>1$ then the
triangular projections are bounded. This implies that the ``tetrahedrals''
are complemented. (An example of this is ${\cal T}=[e_\alpha\colon
\alpha_1<\alpha_2<\cdots<\alpha_N]$). Since all of them are isometrically
isomorphic and there are finitely many of them we conclude that
$X\approx [e_\alpha:\alpha_1<\alpha_2<\cdots<\alpha_n]$ by Pe\l czy\'nski's decomposition method. Then the
proofs of Theorem 3.1 and Proposition 3.2 apply to this space.\endpf
\medbreak
The proof of Theorem 3.1 dualizes (formally) to 
$X_*=\ell_{q_1}\check{\otimes}\cdots\check{\otimes}\ell_{q_N}$.
Define $J_\sigma$ and $K_\sigma$ 
on $B(\ell_{q_1}\check{\otimes}\cdots\check{\otimes}\ell_{q_N})$ 
as in $B(\SN)$. The key to the dualization argument is that
$(J_\sigma)^*=K_\sigma$ and $(K_\sigma)^*=J_\sigma$.

\medbreak
{\sl
THEOREM 3.4. Let $X_*=\ell_{q_1}\check{\otimes}\cdots\check{\otimes}\ell_{q_N}$
be such that ${1\over q_i}+{1\over q_j}\geq1$ for every $i\neq j$. Then
$X_*$ is primary.
}
\medbreak
PROOF. Let 
$\Phi\in B(\ell_{q_1}\check{\otimes}\cdots\check{\otimes}\ell_{q_N})$,
and $0<\epsilon<{1\over 2}$. Then $\Phi^*\in B(\SN)$ and whenever $i\neq j$
we have ${1\over p_i}+{1\over p_j}\leq1$. Therefore, Theorem 3.1 tells us
that there exist $\sigma$ and $\lambda\in{\bf C}$ such that
$\|K_\sigma\Phi^* J_\sigma-\lambda I_X\|<\epsilon$.

Since $K_\sigma\Phi^* J_\sigma-\lambda I_X=
(K_\sigma\Phi J_\sigma-\lambda I_{X_*})^*$ we have that
$\|K_\sigma\Phi J_\sigma-\lambda I_{X_*}\|<\epsilon$. Therefore, $\Phi$
or $I_{X_*}-\Phi$ factors through $X_*$, and since $X_*$ is isomorphic
to its $r'$-sum, we conclude that $X_*$ is primary.\endpf

\bigbreak
\noindent{\bf 4. $\ell_p$ subspaces of $\SN$.}
\medbreak

Theorem 1.3 tells us that $\ell_p$ embeds into $\SN$ if there
exists a non-empty $A\subset\{1,\cdots,N\}$ for which $p=r_A$, where
${1\over r_A}=\min\{1,\sum_{i\in A}{1\over p_i}\}$. We will see in the next
theorem that the converse holds.
\medbreak
{\sl
THEOREM 4.1. $\ell_p$ embeds into $\SN$ if and only if there exists a non-empty
$A\subset\{1,\cdots,N\}$ such that $p=r_A$.
}
\medbreak
We will use this theorem to prove the following:
\medbreak
{\sl
THEOREM 4.2. Let $X=\SN$ and assume that for some $i\neq j$, ${1\over p_i}+
{1\over p_j}>1$ and that $p_k\not\in\{r_A:k\not\in A\}$ for $k=i,j$.
Then $X$ is not primary.
}

\medbreak

REMARK. Theorem 4.1 could probably be generalized to characterize when
$m$-fold tensor products embed into $n$-fold tensor products for $m\leq n$
and this would slightly improve Theorem 4.2.
 
\medbreak
We will use Theorem 1.3 to decompose 
$X\approx [e_\alpha:\alpha_i>\alpha_j]
\oplus [e_\alpha:\alpha_i\leq\alpha_j]$.
The condition $p_i\not\in\{r_A: i\not\in A\}$
insures that $\ell_{p_i}$ does not embed into
$[e_\alpha:\alpha_i\leq\alpha_j]$. (This is easily seen for example when $N=2$.
In this case $X=\ell_{p_1}\hat\otimes\ell_{p_2}$ and 
$p_k\not\in\{r_A:k\not\in A\}$ means
$p_1\neq p_2$; we can then observe that $\ell_{p_1}$ does not
embed into $[e_i\otimes e_j : i\leq j]$).
\medbreak

We will prove Theorem 4.1 by induction. 
Assume for the remainder of this section that $X=\SN$; 
${1\over r}=\min\{1,\sum_{i=1}^N{1\over p_i}\}$ and that
$\Phi:\ell_p\to\SN$ is an isomorphism.
We can assume without loss of generality that there is a sequence of
increasing natural numbers $n_i$ such that
$$\Phi e_i\in[e_\alpha:n_i<|\alpha|<n_{i+1}]
\quad\hbox{ for every  }i.\leqno{(*)}$$
If $p>1$ this is true because $\Phi e_j\to 0$ weakly. If $p=1$ and
$P_{M_n}$ is the projection onto $[e_\alpha:\alpha\leq n]$, we can find
infinitely many pairs of $e_i$'s (say $e_k$ and $e_l$) such that 
$P_{M_n}\Phi(e_k-e_l)\approx 0$. Then we replace the $e_i$'s by differences
of unit vectors and get $(*)$.

We say that $\Psi:\ell_p\to\ell_p$ is an $\ell_p$-average isometry
if there exist a sequence of subsets of $\N$, $\sigma_1<\sigma_2<\cdots$ and
scalars $a_k$ such that
$$\Psi e_i=\sum_{k\in\sigma_i} a_ke_k\quad\hbox{and}\quad
\sum_{k\in\sigma_i}|a_k|^p=1\quad\hbox{for every }i.$$

Finally we will let $E_n=[e_\alpha:\min\{\alpha\}\leq n]$ for every $n\in\N$.
The key to the induction step is that
$$E_n\approx\bigl(\ell_{p_2}\hat\otimes\cdots\hat\otimes\ell_{p_N}\bigr)\oplus
\bigl(\ell_{p_1}\hat\otimes\ell_{p_3}\hat\otimes\cdots\hat\otimes\ell_{p_N}
\bigr)\oplus\cdots\oplus
\bigl(\ell_{p_1}\hat\otimes\cdots\hat\otimes\ell_{p_{N-1}}\bigr).$$
(The isomorphism constant goes to infinity
with $n$). Notice that each one of those summands is an $(N-1)$-projective
tensor product.

We need two lemmas.

\medbreak
{\sl
LEMMA 4.3. Let $\Phi:\ell_p\to\SN$ be as in $(*)$ with
$p>r$. Then for every $\epsilon>0$ we can find $n\in\N$ such that
$\|(I-P_{E_n})\Phi\|\leq\epsilon$.
}
\medbreak

\medbreak
{\sl
LEMMA 4.4. Let $\Phi:\ell_p\to\SN$ be as in $(*)$ with $p<r$, then
for every $\epsilon>0$ there exists $\Psi:\ell_p\to\ell_p$ an $\ell_p$-average
isometry such that $\|\Phi\Psi\|<\epsilon$.
}
\medbreak

PROOF OF THEOREM 4.1.  The theorem is clearly true for $N=1$. Assume that the
result is true for $(N-1)$-projective tensor products and let
$\Phi:\ell_p\to\SN$ be an isomorphism satisfying $(*)$.

It follows from Lemma 4.4 that $p\geq r$. If $p=r$ there is nothing to prove
since $\ell_p$ clearly embeds in the main diagonal. If $p>r$, Lemma 4.3
tells us that $\Phi\ell_p$ is essentially inside $E_n$ and therefore it is
inside one of the $(N-1)$-tensor products. Hence it has to be of the
form $r_A$ for some nonempty $A$ by induction.\endpf
\medbreak
We used in the proof the well-known fact that if $\ell_p$ embeds into $X\oplus Y$ then
$\ell_p$ embeds into $X$ or into $Y$. 
\medbreak
For the proof of Theorem 4.2 we need one more lemma.
\medbreak
{\sl
LEMMA 4.5. Let $X=\SN$, $i,j\leq N$, $i\neq j$ and assume that 
$p_i\not\in\{r_A:i\not\in A\}$. Then
$\ell_{p_i}$ does not embed into $[e_\alpha:\alpha_i\leq\alpha_j]$.
}
\medbreak

PROOF OF THEOREM 4.2.  Use Theorem 1.3 to decompose 
$X\approx [e_\alpha:\alpha_i>\alpha_j]
\oplus [e_\alpha:\alpha_i\leq\alpha_j]$. 
Lemma 4.5
tells us that $\ell_{p_j}$ does not embed
into $[e_\alpha:\alpha_i>\alpha_j]$, and that $\ell_{p_i}$ does not embed
into $[e_\alpha:\alpha_i\leq \alpha_j]$. Therefore neither of them is
isomorphic to $X$, and so $X$ is not primary.\endpf
\medbreak

PROOF OF LEMMA 4.3. If the lemma were false, we could find some
$\epsilon_0>0$; a sequence of normalized vectors $\{x_i\}_{i\in\N}$ 
in $\ell_p$ satisfying
$\hbox{supp}\{x_i\}<\hbox{supp}\{x_{i+1}\}$ for every $i$;
and an increasing sequence $n_i\in\N$ satisfying
$\|P_i\Phi x_i\|\geq\epsilon_0$ where $P_i$ is the projection onto
the diagonal block $[e_\alpha:n_i\leq\alpha<n_{i+1}]$.

Theorem 3.1 implies that $[P_i\Phi x_i:i\in\N]\approx\ell_r$.
Let $P$ be the diagonal projection onto $[P_i\Phi x_i:i\in\N]$
and consider
$P\Phi:\ell_p\to\ell_r$. Since $\|P\Phi e_i\|\geq\epsilon_0$ for
every $i\in\N$ we have that $P\Phi$ is not compact. This 
is a contradiction.\endpf
\medbreak

SKETCH OF THE PROOF OF LEMMA 4.4. For $N=1$ the result is easy. 
The condition $(*)$ says that 
$\Phi:\ell_p\to\ell_r$ is diagonal; i.e.,
$\Phi e_i=\lambda_i e_i$. Moreover since $\Phi$ is bounded, there exists $M>0$
such that 
$|\lambda_i|\leq M$ for every $i$. 
We get the blocks by taking the $a_k$'s constant in every $\sigma$. Let
$\sigma\subset\N$ be of cardinality $n$ (say). Then
$\|\sum_{k\in\sigma}\bigl({1\over n}\bigr)^{1/p} e_k\|_p=1$ but
$\|\sum_{k\in\sigma}\bigl({1\over n}\bigr)^{1/p}\Phi e_k\|_r\leq M n^{1/r-1/p}$
goes to
zero as $n$ goes to infinity.

Assume the result for $N-1$ and let $\Phi:\ell_p\to\SN$ be as in
$(*)$. The idea is to find an $\ell_p$-average isometry $\Psi\in B(\ell_p)$
such that $\Phi\Psi$ is essentially supported in a diagonal block; then
since the diagonal block is like $\ell_r$, the case $N=1$ takes care of it.

To find $\Psi$ we have to find an increasing sequence $n_i\in\N$ and
a normalized sequence $\{x_i\}_{i\in\N}$ in $\ell_p$ satisfying
supp$\{x_i\}\leq\hbox{ supp}\{x_{i+1}\}$ for every $i\in\N$ and 
$\Phi x_i\in [e_\alpha:n_i\leq\alpha<n_{i+1}]$. (The last inclusion is an
``almost'' inclusion; that is, for a given $\epsilon_i>0$ 
there exists $n_i\in\N$ such that the distance from
$\Phi x_i$ to $[e_\alpha:n_i\leq\alpha<n_{i+1}]$ is less that $\epsilon_i$).

It is clear that it is enough to do this for $x_1$ and $x_2$ because
we can iterate it to conclude the lemma. Clearly 
$\Phi x_1\in [e_\alpha:
\alpha<n]$ for some $n$. We want to find $x_2$ such that
$\Phi x_2$ is supported outside $E_n$. Since $E_n$ is isomorphic to the
sum of $(N-1)$-projective tensor products, we can apply the induction
step to insure the existence of $x_2$. \endpf
\medbreak
SKETCH OF THE PROOF OF LEMMA 4.5. The proof of this goes by induction too.
The result is clear for $N=2$.  Suppose it is true for $N-1$ and 
false for $N$. Then let $Z=[e_\alpha:\alpha_i\leq
\alpha_j]\subset\SN$ and, by the assumption, find $\Phi:\ell_{p_i}\to Z$, an isomorphism satisfying $(*)$. 

The main diagonal of $Z$ is isomorphic to $\ell_r$ and $p_i>r$.
Hence, by Lemma 4.3, there exists $n\in\N$ such that $\Phi\ell_{p_i}$ is essentially
inside $E_n$. We will look at the $N$-summands of $Z\bigcap E_n$ to get a
contradiction.

One of those summands does not contain the $i$th component and hence 
is isomorphic
to $\ell_{p_1}\hat\otimes\cdots\hat\otimes\ell_{p_{i-1}}\hat\otimes
\ell_{p_{i+1}}\hat\otimes\cdots\hat\otimes\ell_{p_N}$. The condition  $p_i\not\in\{r_A: i\not\in A\}$ and Theorem 4.1 imply that
$\ell_{p_i}$ does not embed there.

Another summand does not contain the $j$th component. This really means
that $\alpha_j\leq n$. Therefore, $\alpha_i\leq n$ as well and the
summand is isomorphic
to $\hat\otimes_{k\neq i,j}\ell_{p_k}$. We conclude as before that
$\ell_{p_i}$ does embed here.

The remaining summands will have the same structure but with $N-1$ terms.
Then the induction hypothesis implies that $\ell_{p_i}$ does not embed
into any one of them.

Therefore, $\ell_{p_i}$ does not embed in $Z$. This is a contradiction.\endpf
\bigbreak

\bigbreak
\noindent{\bf 5. Primarity of Polynomials and Operator Spaces.}
\medbreak
In this section we discuss the primarity of $(\SN)^*$. There
will be really only one case to consider; namely that of $r=1$
(recall that ${1\over r}=\min\{1,\sum_{i=1}^N{1\over p_i}\})$,
which we demonstrate below using techniques of Bourgain [Bo] and Blower [B].

It is interesting to note that completely different factors determine
the primarity of $(\SN)^*$ when $r=1$ and $r>1$. When $r>1$ it is the unboundedness of the main triangle projection in each pair
(taken separately) that is the most important factor, while for $r=1$
we will see that the main point is that  
we have $\ell_\infty$-blocks down the diagonal.
\medbreak
{\sl
THEOREM 5.1. Let $X=\SN$ be such that 
${1\over r}=\min\{1,\sum_{i=1}^N{1\over p_i}\}=1$. Then $(\SN)^*$ is primary.
}
\medbreak
This result will solve the question of primarity for spaces of polynomials.
Since the space of analytic polynomials of degree $m$ on $\ell_p$ is isomorphic 
(with constant ${m^m\over m!}$) to the dual of the symmetric $m$-fold tensor
product ${\hat{\otimes}}_s^m\ell_p$. That is ${\cal P}_m\approx
({\hat{\otimes}}_s^m\ell_p)^*$.  (Here $m$ is the
number of times that one takes the tensor product).
\medbreak
{\sl
LEMMA 5.2. For any $1\leq p<\infty$ and $m\in\N$ we have that
$\ell_p\hat\otimes\cdots\hat\otimes\ell_p\approx{\hat{\otimes}}_s^m\ell_p$.
}
\medbreak
PROOF. We use Pe\l czy\'nski's decomposition method again. Since
$\ell_p\hat\otimes\cdots\hat\otimes\ell_p$ is isomorphic to its infinite
$s$-sum ($s=\max\{1,{p\over m}\}$) we only have to prove that they embed 
complementably into each other. It is clear that ${\hat{\otimes}}_s^m\ell_p$
embeds into $\ell_p\hat\otimes\cdots\hat\otimes\ell_p$. Indeed,
$S\in B(\ell_p\hat\otimes\cdots\hat\otimes\ell_p)$ defined by
$Se_\alpha={1\over m!}\sum_{\pi\in \Pi_m}e_{\pi(\alpha)}$ 
shows that the embedding is 1-complemented. On the other hand,  
for $i\leq N$ let $\sigma_i(j)=m(j-1)+i$, 
$\sigma=(\sigma_1,\cdots,\sigma_m)$ and define 
$T\in B(\ell_p\hat\otimes\cdots\hat\otimes\ell_p)$ by $T=K_\sigma SJ_\sigma$.
It is clear that $T$ factors through ${\hat\otimes}_s^m\ell_p$ and 
it is easy to see that $Te_\alpha={1\over m!}e_\alpha$. Hence, 
$\ell_p\hat\otimes\cdots\hat\otimes\ell_p$ embeds complementably into
${\hat\otimes}_s^m\ell_p$ and the result follows.\endpf
\medbreak
{\sl
COROLLARY 5.3. Let $1\leq p<\infty$ and $m\geq1$. The space of homogeneous
analytic polynomials ${\cal P}_m(\ell_p)$ and the symmetric tensor
product of $m$ copies of $\ell_p$ are primary.
}
\medbreak
We now proceed to the proof of the theorem. Notice that if
$X=\SN$ is such that 
${1\over r}=\min\{1,\sum_{i=1}^N{1\over p_i}\}=1$, 
then Theorem 2.2 tells us that
$$(\SN)^*\approx \bigl(\sum_{n=1}^\infty
\ell_{q_1}^n\check{\otimes}\cdots
\check{\otimes}\ell_{q_N}^n\bigr)_\infty.$$

This decomposition allows us to use the technique developed by Bourgain [Bo]
to prove that $H^\infty$ is primary; namely, one obtains the general
theorem from the finite dimensional version.

The proof is an exact
generalization of the proof of Blower [B] that $B(H)$ is
primary; it has no surprises, and so we will simply sketch
the part that is different for the case $N>2$, and refer
the interested reader to [Bl
] for other details. The proof follows
from the following 2 lemmas, as indicated in [Bo].
\medbreak
{\sl
PROPOSITION 5.4. Given $n\in{\bf N}$, $\epsilon>0$ and $K<\infty$, there exists $N_0=N_0(n,\epsilon,K)$ such that if $M\geq N_0$ and  
$T\in B(\ell_{q_1}^M\check{\otimes}\cdots\check{\otimes}\ell_{q_N}^M)$
with $\|T\|\leq K$, then there exist subsets 
$\sigma_1,\sigma_2,\cdots\sigma_N\subset\{1,\cdots,M\}$ of cardinality $n$, 
and a constant $\lambda$ such that if $\sigma=\Nsigma$ then, 
$$\left\|K_\sigma TJ_\sigma-\lambda I_n\right\|\leq\epsilon.$$
Thus, one of $K_\sigma TJ_\sigma$ and 
$K_\sigma(I_N-T)J_\sigma$ is invertible.
}
\medbreak
REMARK. Here $J_\sigma\colon \ell_{q_1}^n\check\otimes\cdots\check\otimes
\ell_{q_N}^n\to\ell_{q_1}^M\check\otimes\cdots\check\otimes\ell_{q_N}^M$ is
defined by $J_\sigma e_\alpha=e_{\sigma(\alpha)}$ where 
$\sigma(\alpha)=(\sigma_1(\alpha_1),\cdots,\sigma_N(\alpha_N))$, and
$\sigma_i=\{\sigma_i(1),\sigma_i(2),\cdots,\sigma_i(n)\}$. Moreover, 
$\sigma_i(k)<\sigma_i(l)$ iff $k<l$. The definition for 
$K_\sigma$ is similar.

\medbreak
{\sl
PROPOSITION 5.5.  Given $n\in{\bf N}$ and $\epsilon>0$ there exists 
$N_0=N_0(n,\epsilon)$ such that if $M\geq N_0$ and $E$ is an 
$n$-dimensional subspace of 
$\ell_{q_1}^M\check{\otimes}\cdots\check{\otimes}\ell_{q_N}^M$ 
then there exists a subspace $F$ of
$\ell_{q_1}^M\check{\otimes}\cdots\check{\otimes}\ell_{q_N}^M$, 
isometrically isomorphic to 
$\ell_{q_1}^n\check{\otimes}\cdots\check{\otimes}\ell_{q_N}^n$, 
and a block projection 
$Q$ from $\ell_{q_1}^M\check{\otimes}\cdots\check{\otimes}\ell_{q_N}^M$ 
to $F$ such that $\|Qx\|<\epsilon\|x\|$ for every $x\in E$.
}
\medbreak
SKETCH OF THE PROOF OF PROPOSITION 5.4. 
Let $T\in B(\ell_{q_1}^M\check\otimes\cdots
\check\otimes\ell_{q_N}^M)$ such that $\|T\|\leq K$. We will find
a copy of $\ell_{q_1}^n\check\otimes\cdots\check\otimes\ell_{q_N}^n$
inside $\ell_{q_1}^M\check\otimes\cdots\check\otimes\ell_{q_N}^M$ such
that $T$ is essentially a multiple of the identity when restricted to this subspace.
We accomplish this in two steps.

STEP 1. Find a large subset $\psi\subset\{1,\cdots,M\}$ and
$\lambda\in{\bf C}$ such that whenever $\alpha=(\alpha_1,\cdots,\alpha_N)$
is such that $\alpha_1<\cdots<\alpha_N$ and 
$\alpha_k\in\psi$ for $i\leq N$, then 
$|(Te_\alpha,e_\alpha)-\lambda|<\epsilon$.

STEP 2. Find $\sigma_1<\sigma_2<\cdots<\sigma_N\subset\psi$ each of
cardinality $n$, such that whenever 
$\alpha=(\alpha_1,\cdots,\alpha_N), 
\alpha'=(\alpha'_1,\cdots,\alpha'_N)$ are 
are such that $\alpha_k,\alpha'_k\in\sigma_k$ for every $k\leq N$ and
$\alpha\neq\alpha'$,
then $|(Te_\alpha,e_{\alpha'})|<\epsilon$. Then define
${\cal S}=[e_\alpha\colon \alpha_k\in\sigma_k]$.  One can easily verify
that if $\epsilon>0$ is chosen small enough then $T$ restricted and projected
into ${\cal S}$ is essentially a multiple of the identity and that
${\cal S}$ is isometrically isomorphic to
$\ell_{q_1}^n\check\otimes\cdots\check\otimes\ell_{q_N}^n$.

Both steps depend on Ramsey's Theorem and they are very
minor modifications of Blower's argument.

For Step 1 divide the disk $\{z\colon|z|\leq K\}$ into finitely many
disjoint subsets $V_k$ of diameter less than $\epsilon$, and define the 
coloring on $N$-sets of $\{1,\cdots,M\}$ by 
$\{\alpha_1,\cdots,\alpha_N\}\to\ell$
if $(Te_\alpha,e_\alpha)\in V_{\ell}$ where
$\alpha=(\alpha_1,\cdots,\alpha_N)$ for $\alpha_1<\cdots<\alpha_N$.
Then use Ramsey's Theorem fo find a large monochromatic set $\psi$.

The proof of Step 2 involves many different cases (but all of them are
similar). One has to look at all the different ways that
$(\alpha_1,\cdots,\alpha_N)\neq(\alpha'_1,\cdots,\alpha'_N)$. We will illustrate the case
when $\alpha_k<\alpha'_k$ for every $k\leq N$.

Color the $2N$-elements of $\{1,\cdots,M\}$ by: $\{\alpha_1,\alpha'_1,\alpha_2,\alpha'_2,\cdots,
\alpha_N,\alpha'_N\}$ is {\it bad} if $\alpha_1<\alpha'_1<\alpha_2<\alpha'_2<\cdots<\alpha_N<\alpha'_N$ and
$|(Te_\alpha,e_{\alpha'})|\geq\epsilon$ where $\alpha=(\alpha_1,\cdots,\alpha_N)$ and
$\alpha'=(\alpha'_1,\cdots,\alpha'_N)$; it is {\it good} otherwise.

Ramsey's Theorem gives us a large monochromatic subset $\psi_1\subset\psi$.
We will show that $\psi_1$ has to be {\it good}. Let $\alpha'_1<\alpha_2<\alpha'_2<\cdots
<\alpha_N<\alpha'_N$ be the $2N-1$ largest elements of $\psi_1$, and let
$\beta=(\alpha_2,\cdots,\alpha_N)$, $\alpha'=(\alpha'_1,\cdots,\alpha'_N)$, and 
$F=[e_i\otimes e_\beta\colon i\in\psi_1, i<\alpha'_1]$.

It is clear that $F\equiv\ell_{p_1}^{|\psi_1|-2N+1}$. Define 
$\tilde{T}:F\to{\bf C}$ by $\tilde{T}(x)=(Tx,e_{\alpha'})$. Then we have that
$\tilde{T}$ is a map from $\ell_{p_1}^s$ into ${\bf C}$, with norm less
than or equal to $K$ and maps the canonical basis into ``large'' elements.
Since we assumed that $p_1>1$ this is a contradiction.

Now we have to look at all the other possibilities; e.g., $\alpha'_1>\alpha_1$ and
$\alpha_k<\alpha'_k$ for $2\leq k\leq N$ etc. We have to look also at the cases when 
some of the coordinates are equal, but these are not very different.
We prove the Proposition by choosing $M$ large enough.\endpf
\medbreak

SKETCH OF THE PROOF OF PROPOSITION 5.5.  It is enough to prove that if
$x\in\ell_{q_1}^M\check\otimes\cdots\check\otimes\ell_{q_N}^M$ 
with $\|x\|\leq1$,
then we can find $Q$, a large block projection, such that $\|Q(x)\|\leq
\epsilon$. Then take an $\epsilon$-net of the sphere of $E$, 
$\{x_i\}_{i=1}^s$. Find $Q_1$ a large block projection such that
$\|Q_1x_1\|\leq\epsilon$; then find $Q_2$ a large block  projection
contained in the range of $Q_1$ such that $\|Q_2Q_1x_2\|\leq\epsilon$.
Proceeding in this way we get that $Q=Q_s\cdots Q_2Q_1$; this $Q$
does it.

To check the first claim let $x\in\ell_{q_1}^M\check\otimes\cdots\check\otimes
\ell_{q_N}^M$ with $\|x\|\leq1$ and let $\rho>0$ (to be fixed later). 
Then define
a coloring on the $N$ sets of $\{1,\cdots,M\}$ by: $\{\alpha_1,\cdots,\alpha_N\}$ is
{\it bad} if $|(x,e_\alpha)|\geq\rho$ where $\alpha=(\alpha_1,\cdots,\alpha_N)$ and
$\alpha_1<\cdots<\alpha_N$. And {\it good} otherwise. Ramsey's theorem gives us a
large monochromatic subset, and this subset has to be {\it good}.\endpf
\bigbreak

\noindent{\bf 6. Appendix}
\medbreak
In this section we will prove Basic Lemmas 1 and 
2 from Section 3. 
\medbreak
PROOF OF BASIC LEMMA 1.  For this proof let $M_n=[e_\alpha:\alpha\leq n]$
with projection $P_n$. We will divide the proof into two parts, one for
$m>n$ and the other one for $m<n$.  In both cases,
$\sigma=(\sigma_1,\cdots,\sigma_N)$ satisfies $\sigma_1=\sigma_2=\cdots
=\sigma_N.$ 
\medbreak
The case $m>n$ is simpler; we start with it.

If $K\subset\SN$ is a compact set, then $K$ is essentially inside one
of the $M_n$'s. The following
elementary lemma states this fact quantitatively (we omit its proof as it is an easy exercise). The proof of the case
$m>n$ follows easily from it.

\medbreak
{\sl
LEMMA 6.1. Let $K\subset\SN$ be a compact set and $\epsilon_k>0$ be given.
Then we can find a sequence $n_k\in{\bf N}$
such that $\sup_{x\in K}\|(I-P_{n_k})x\|<\epsilon_k$.
}
\medbreak

We start the inductive construction of $\sigma_1$. Set $A_1={\bf N}$ and 
$\sigma_1(1)=\min A_1$. Let $K=\Phi \hbox{ Ball }\partial M_{\sigma_1(1)}$
and $\epsilon_k=\epsilon_{1,k}$. Then find 
$A_2\subset A_1\setminus\{\sigma_1(1)\}$ according to Lemma 6.1;
and set $\sigma_1(2)=\min A_2$.

Let $K=\Phi \hbox{ Ball }\partial M_{\sigma_1(2)}$, 
$\epsilon_k=\epsilon_{2,k}$ and find
$A_3\subset A_2\setminus\{\sigma_1(2)\}$ according to Lemma 6.1. Then set
$\sigma_1(3)=\min A_3$.

Continuing in this fashion we get $\sigma_1$ and construct $\sigma=(\sigma_1,
\cdots,\sigma_1)$. It is easy to see that if $x\in\partial M_n$ and $m>n$, then
$$\|Q_mK_\sigma\Phi J_\sigma x\|\leq\epsilon_{n,m}\|x\|.$$
\medbreak 

We will now prove the case $m<n$. 

The construction of $\sigma_1$ is similar to the previous case. We need     
the following elementary lemma (which as before do not prove).
\medbreak
{\sl
LEMMA 6.2. Let $1<p<\infty$, $F$ a finite dimensional space,
and $T:\ell_p\to F$ a bounded linear map. Then for every $\epsilon>0$, the
set $\{i:\|Te_i\|>\epsilon\}$ is finite.
}
\medbreak

We will only present the induction step for the construction of $\sigma_1$.
Assume that $\Lambda\subset{\bf N}$ is an infinite set with first $n$ 
elements $\sigma_1(1),\cdots,\sigma_1(n)$. We want to find  an infinite
$\Lambda'\subset\Lambda$ with the same first $n$ elements as $\Lambda$
such that whenever $\alpha\in(\Lambda')^N$ is such that
$e_\alpha\not\in M_{\sigma_1(n)}$, then $\|P_{\sigma_1(n)}\Phi e_\alpha\|<
\epsilon$. Then we will choose $\sigma_1(n+1)=\min \Lambda'\setminus
\{\sigma_1(1),\cdots,\sigma_1(n)\}$.

The construction of $\Lambda'$ uses Ramsey's Theorem as in Section 5. We look
at all the different ways that $e_\alpha\not\in M_{\sigma_1(n)}$. We will
illustrate this for two different cases. The others are very similar.

\smallbreak
CASE 1: $\sigma_1(n)<\alpha_1<\alpha_2<\cdots<\alpha_N$.

Color the $N$-sets of $\{i\in\Lambda:i>\sigma_1(n)\}$ as follows: 
$\{\alpha_1,\cdots,\alpha_N\}$ is {\it good} if $\alpha_1<\cdots<\alpha_N$
and $\|P_{\sigma_1(n)}\Phi e_\alpha\|<\epsilon$ and {\it bad} otherwise.

Ramsey's Theorem gives us a monochromatic 
infinite set $\Lambda_1\subset\Lambda$. It is easy to see that Lemma 6.2
implies that the set has to be {\it good}. (Let $\beta_1<\cdots
<\beta_{N-1}$ be the $N-1$ smallest elements of $\Lambda_1$ and define
$T:\ell_{p_N}\to M_{\sigma_1(n)}$ as follows: if $i>\beta_{N-1}$, then 
$T e_i=P_{\sigma_1(n)}\Phi e_{(\beta_1,\cdots,\beta_{N-1},i)}$ and
if $i\leq\beta_{N-1}$, then $Te_i=0$.
If $\Lambda_1$ were {\it bad} this would contradict Lemma 6.2.)
\smallbreak
CASE 2: $\alpha_1,\alpha_2\leq\sigma_1(n)<\alpha_3<\cdots<\alpha_N$.

Color the $(N-2)$-sets of $\{i\in\Lambda:i>\sigma_1(n)\}$ as follows:
$\{\alpha_3,\cdots,\alpha_N\}$ is {\it good} if $\alpha_3<\cdots<\alpha_N$
and $\|P_{\sigma_1(n)}\Phi e_\alpha\|<\epsilon$ for every $\alpha_1,\alpha_2
\leq \sigma_1(n)$. (Notice that $\alpha=(\alpha_1,\alpha_2,\alpha_3,\cdots,
\alpha_N)$) and {\it bad} otherwise.

Once again Ramsey's Theorem gives an infinite monochromatic subset of $\Lambda$.
And as before it has no choice but to be {\it good}. This follows because there are
only finitely many $\alpha_1,\alpha_2\leq\sigma_1(n)$.
\medbreak
There are finitely many ways in which $e_\alpha\not\in M_{\sigma_1(n)}$.
They are very similar to the two cases just considered, and repeating the above argument for 
all of them we get $\tilde{\Lambda}\subset\Lambda$ that is {\it good} in
all the cases. Then let $\Lambda'=\tilde{\Lambda}\bigcup
\{\sigma_1(1),\cdots,\sigma_1(n)\}$.
We choose $\epsilon>0$ small enough so that whenever $x\in\partial M_{n+1}$,
then 
$$\|P_{\sigma_1(n)} K_\sigma \Phi J_\sigma x\|
\leq\min_{k\leq n}\epsilon_{n+1,k}\|x\|.\endpf$$

\bigbreak
PROOF OF BASIC LEMMA 2. Assume that we have a sequence of complex numbers
$\{\lambda_{\alpha,\beta,j}:\alpha,\beta\in{\bf N}^{N-1}, |\alpha|\vee|\beta|<j\}$
and a sequence of positive numbers,
$\{\epsilon_{\alpha,\beta,j}:j>|\alpha|\vee|\beta|\}$.

For $\alpha,\beta$ fixed, find a subsequence $\{j_k\}$
of $\{j:j>|\alpha|\vee|\beta|\}$ and some $\lambda_{\alpha,\beta}\in{\bf C}$ 
satisfying:
$$\eqalign{\lim_{j_k\to\infty}\lambda_{\alpha,\beta,j_k}
                                            &=\lambda_{\alpha,\beta}\cr
           |\lambda_{\alpha,\beta,j_k}-\lambda_{\alpha,\beta}|
                          &<\epsilon_{\alpha,\beta,k}.\cr}\leqno{(**)}$$
Moreover, if we have finitely many 
$\{\alpha_l,\beta_l\}_{l\leq m}$, we can find a subsequence $\{j_k\}$
such that $(**)$ is true for every $l\leq m$.

The condition $j>|\alpha|\vee|\beta|$ is the key to extend the argument to
all $\alpha,\beta\in{\bf N}^{N-1}$. The basic idea is that once we have fixed
$\sigma_1(1),\cdots,\sigma_1(n)$, 
we take the subsequence $j_k$ from $\{j:j>\sigma_1(n)\}$; hence, we
do not affect the initial segment. 

We will only present the induction step for $\sigma_1$. Assume that
$\Lambda\subset{\bf N}$ is an infinite set with first elements 
$\sigma_1(1),\sigma_1(2),\cdots,\sigma_1(n)$. We want to find an infinite
$\Lambda'\subset\Lambda$ with the first $n$ elements as in $\Lambda$, and such
that
$(**)$ is satisfied for every $\alpha,\beta\leq\sigma_1(n)$. We can do that
because there are only finitely many of them. We take the subsequence $j_k$
from $\{j\in\Lambda:j>\sigma_1(n)\}$ and let 
$\Lambda'=\{j_k:k\in{\bf N}\}\bigcup\{\sigma_1(1),\cdots,\sigma_1(n)\}$.
Then set $\sigma_1(n+1)=j_1$, the minimum of the $j_k$'s (remember that
$j_1>\sigma_1(n)$).

Repeating the process we finish the proof.\endpf
\bigbreak
\centerline{REFERENCES}

\medbreak
\item{[AEO]} D. Alspach, P. Enflo and E. Odell, {\it On the structure of ${\cal L}_p$, $(1<p<\infty)$}, Studia Math. {\bf 60} (1977), 79--90.

\medbreak
\item{[Ar1]}
J. Arazy, {\it Basic sequences, embeddings, and the uniqueness
of the symmetric structure in unitary spaces}, J. Func. Anal. {\bf  
40} (1981),
302--340.

\medbreak
\item{[Ar2]}
J. Arazy, {\it On subspaces of $c_p$ which contain $c_p$},
Compositio Math. {\bf 41} (1980), 297--336.

\medbreak
\item{[A]} A. Arias, {\it Nest algebras in the trace class},  
Journal of Operator Theory (to appear).

\medbreak
\item{[Be]} G. Bennett, {\it Unconditional convergence and almost  
everywhere convergence}, Z. Wah. verw. Gebiete  
{\bf 34} (1976) 135--155.

\medbreak
\item{[B]} G. Blower, {\it The space $B(H)$ is primary}, Bull.  
London Math. Soc. {\bf 22} (1990), 176--182.

\medbreak
\item{[Bo]} J. Bourgain, {\it On the primarity of ${\cal  
H}^\infty$-spaces}, Israel J. Math {\bf 45} (1983), 329--336.

\medbreak
\item{[CL]} P. Casazza and B. Lin, {\it Projections on Banach spaces with symmetric bases}, Studia Math. {\bf 52} (1974), 189--193.

\medbreak
\item{[FJ]} J. Farmer and W. B. Johnson, {\it Polynomial Schur and polynomial Dunford-Pettis properties}, Banach Spaces, Contemporary Mathematics. {\bf 144}  (1993), Bor-Luh Lin, (ed.), American Mathematical Society,  Providence, RI.

\medbreak
\item{[JRZ]} W. B. Johnson, H. P. Rosenthal, and M. Zippin, {\it On
bases, finite dimensional decompositions, and weaker structures in Banach
spaces}, Israel J. Math. {\bf 9} (1971), 488--506.

\medbreak
\item{[KP]} S. Kwapien and A. Pe\l czy\'nski, {\it The main  
triangular
projection in matrix spaces and its applications}, Studia Math. {\bf 34} (1970),
43--68.

\medbreak
\item{[LP]} J. Lindenstrauss and A. Pe\l czy\'nski,  {\it Absolutely summing operators in 
${\cal L}_p$--spaces and their applications}, Studia Math. {\bf 29} (1968) 
{275--326}.

\medbreak
\item{[LT]} J. Lindenstrauss and L. Tzafriri, {\it Classical Banach
spaces I: Sequence spaces},  Springer-Verlag, {Berlin} (1977).   

\medbreak
\item{[MN]} B. Maurey and A. Nahoum: Note aux C.R. Acad. Sc. Paris, t. 276, March 1973, page 751.

\medbreak
\item{[M]} J. Mujica, {\it Complex analysis in Banach spaces,} Notas de  
Mathematica, vol. 120, North-Holland, Amsterdam, 1986.

\medbreak
\item{[R]}
R. Ryan, {\it Applications of topological tensor products to infinite  
dimensional holomorphy}, Ph.D. Thesis, Trinity College, Dublin, 1980.

\medbreak
\item{[Z]} I. Zalduendo, {\it An estimate for multilinear  
forms on $\ell_p$ spaces}, preprint.

\bye